\documentclass{article}
\usepackage{amsmath,amsfonts,amssymb}
\usepackage{graphicx}
\usepackage{eepic}
\usepackage{epic}

\newtheorem{theorem}{Theorem}[section]
\newtheorem{definition}[theorem]{Definition}

\newtheorem{example}[theorem]{Example}
\newtheorem{lemma}[theorem]{Lemma}

\newtheorem{remark}[theorem]{Remark}

\newcommand{\mE}{\mathbb E}

\newcommand{\mN}{\mathbb N}

\newcommand{\mC}{\mathbb C}

\newcommand{\lC}{\cal C}

\newcommand{\mR}{\mathbb R}

\newcommand{\be}{\begin{eqnarray}}
\newcommand{\ee}{\end{eqnarray}}
\newcommand{\bd}{\begin{definition}}
\newcommand{\ed}{\end{definition}}

\newcommand{\br}{\begin{remark}}
\newcommand{\er}{\end{remark}}
\newcommand{\gog}{{\mathfrak g}}
\newcommand{\gob}{{\mathfrak b}}

\newcommand{\gop}{{\mathfrak p}}

\newcommand{\gon}{{\mathfrak n}}
\newcommand{\gol}{{\mathfrak l}}

\newcommand{\bt}{\begin{tabular}}
\newcommand{\et}{\end{tabular}}

\input xypic
\input epsf
\input xy
\xyoption{all}

\begin{document}
\baselineskip13pt
\date{}
\title{Branching problems and ${\mathfrak sl}(2,\mC)$-actions}
\author{ Pavle Pand\v zi\' c, Petr Somberg
}

\maketitle

\abstract 

We study certain ${\mathfrak sl}(2,\mC)$-actions associated to specific examples of
branching of scalar generalized Verma modules for compatible pairs $(\gog,\gop)$, 
$(\gog',\gop')$ of Lie algebras and their parabolic subalgebras. 

\vspace{0.3cm}

{\bf Key words:} Representation theory of simple Lie algebra, Generalized
Verma modules, Singular vectors and composition series, Relative Lie algebra 
and Dirac cohomology.
%
\endabstract


\section{Introduction}

The notion of composition series or branching rules for 
geometrically realized algebraic objects of representation 
theoretical origin, lies at the heart of many 
problems on the intersection of representation theory, algebraic analysis 
and differential geometry.  

The present letter attempts to address several concrete questions
belonging to this line of research, originating in the series 
of articles \cite{koss}, \cite{koss2}, \cite{kp}. Namely, motivated
by questions in differential geometry and harmonic analysis on 
differential invariants associated to pairs of generalized flag manifolds,
a constructive method (the F-method) was developed there. It is based 
on algebraic analysis applied to generalized Verma modules, realized 
as ${\fam2 D}$-modules on the point orbit of the nilpotent group on 
the generalized flag manifold. The methods of algebraic analysis,
which replace the standard combinatorial approach, allow to find the
singular vectors responsible for the composition structure of generalized 
Verma modules in a striking way. In many cases, the factorization 
identities yield the answer 
not only in the Grothendieck group of the Bernstein-Gelfand-Gelfand 
category ${\fam2 O}^\gop(\gog)$, but also allow to recognize the 
composition structure and the extension 
classes.  

In the present note we discuss several questions left untouched in \cite{koss},
\cite{koss2}. To describe these questions, we first introduce some notation. 
Let $\gog$ be a simple Lie algebra and let $\gop\subset\gog$ be a parabolic subalgebra.
Let  $(\gog', \gop'), (\gog, \gop)$ be a compatible pair of Lie algebras and 
their parabolic subalgebras with 
$\gog'\subset\gog, \gop'\subset\gop$.
Let $\gon$, respectively $\gon'$, be the nilradical of $\gop$ and $\gop'$,
respectively.

We focus on the role of the generators of the 
complement of $\gon'$ in $\gon$.
As we shall see in our two examples, each of them is characterized by the 
one dimensional quotient $\gon/\gon'$, the root vector generating 
this complement acts on the space of $\gog'$-singular vectors. We 
shall observe in one of these examples that there is moreover a ${\mathfrak sl}(2,\mC)$-module
structure on the space of $\gog'$-singular vectors, but the action of the opposite 
root space does not directly follow from any representation theoretical construction. 

The article is organized as follows. In the beginning of Section $2$
we review some basic notation, and we pass to the two examples -- one is given
by the Hermitean symmetric space associated to the first node of the Dynkin diagram
of an orthogonal Lie algebra, and the second example is the Lie algebra
${\mathfrak sl}(2,\mR)$ diagonally embedded in 
${\mathfrak sl}(2,\mR)\oplus{\mathfrak sl}(2,\mR)$. As we already mentioned,
in the first example we find an ${\mathfrak sl}(2,\mC)$-module structure on 
the set of $\gog'$-singular vectors which allows detailed analysis related 
to the structure of composition series.
In Section $3$, we pass to a natural question on the relative Lie algebra or Dirac 
cohomology associated to our branching problem, which does not seem to be discussed
in the literature. We finish the letter by several useful conventions and 
formulas related to Gegenbauer and Jacobi polynomials, which realize the 
$\gog'$-singular vectors in the examples of Section $2$.


\section{Main examples}

We shall start with a brief review of several basic notions, 
relying on the conventions 
in \cite{koss}, \cite{koss2}. 

We denote by $\gog$, $\gop$, $\gol$, $\gon$, $\gon_-$ a real simple Lie algebra, its 
parabolic subalgebra, the Levi factor and the nilradical of the parabolic subalgebra, and 
the opposite (negative) nilradical. There is an isomorphism of 
vector spaces 
$\gog\simeq \gon_-\oplus\gop\simeq \gon_-\oplus\gol\oplus\gon$.
The connected and simply connected groups corresponding to
these Lie algebras are denoted by $G,P,L,N,N_-$.

The $\gog$-modules we consider are the scalar generalized Verma modules, 
geometrically realized by ${\fam2 D}$-modules
supported at the closed orbit $eP$ of the nilpotent group $N$ on 
the generalized flag manifold $G/P$. 
These modules are then identified by distribution Fourier transform 
with the 
underlying vector space of the polynomial algebras in Fourier dual
variables $\xi_i$, $i=1,\cdots , \mathrm{dim}_\mR(\gon_-)$.
We can consider another collection of Lie algebras 
$\gog'$, $\gop'$, $\gol'$, $\gon'$, $\gon_-'$
with the same properties as in the un-primed case, such that 
$\gog'\subset\gog$ induces the 
primed to un-primed inclusions for all the other Lie subalgebras.  
The compatibility of the two collections through
the compatible grading of $\gop'$ and $\gop$ realized by the 
adjoint action of the grading element in $\gol'$ is 
also required. 

The $\gog'$-singular vectors in $M^\gog_\gop(\lambda)$ 
describe the generators of $\gog'$-submodules in the 
Grothendieck group $K({\fam2 O}^{\gop'})$ of the 
Bernstein-Gelfand-Gelfand parabolic category ${\fam2 O}^{\gop'}$,
and consequently determine its $\gog'$-composition structure. 
They are the solution
spaces of the system of partial differential equations corresponding to
the action of $\gon'$ by differential operators on 
$\mC[\xi_1,\ldots , \xi_{\mathrm{dim}_\mR(\gon_-)}]$. 

In the two examples of our interest, the $\gog'$-singular vectors 
are the Gegenbauer and Jacobi polynomials, respectively. Our main
concern in these simple (but representative) examples is the action 
of the generator of $\gon/\gon'$ on the span 
of $\gog'$-singular vectors and its consequences, e.g. a lift of its 
structure to a 
${\mathfrak sl}(2,\mC)$-structure. We construct such a lift in one of 
our examples, but at the cost of having to leave the setting of the universal enveloping algebra $U(\gog)$. 

This raises a representation theoretical question whether, say even
for the Hermitean symmetric spaces $(\gog,\gop)$ characterized by parabolic 
subalgebras $\gop$ with commutative nilradical, the action of elements 
in $\gon/\gon'$ on the $U(\gol')$-submodule $\mathrm{Ker}(\gon')\subset M^\gog_\gop(\lambda)$ of a 
generalized Verma module can be lifted to a non-trivial action of a bigger 
subalgebra in the Weyl algebra on the opposite nilradical $\gon$. 
We do not know the answer
to this question, and the present article is a modest attempt to get an
elementary insight into some of its aspects.

Throughout the article we use the notation $\langle \, ,\rangle$ for the linear span of a 
subset of a vector space, $U$ applied to an algebra denotes its universal 
enveloping algebra, and $\mN_0$ the set of natural integers including $0$.

\subsection{The pair of Lie algebras ($\mathfrak{so}(n+1,1,\mR), \mathfrak{so}(n,1,\mR)$) 
and the Gegenbauer polynomials}
\label{subsection2.1}
Let us consider the case of compatible pair of Lie algebras 
\begin{eqnarray}
& & \gog=\mathfrak{so}(n+1,1,\mR),\, \gop=(\mathfrak{so}(n,\mR)\times\mR)\ltimes\mR^n,
\nonumber \\
& & \gog'=\mathfrak{so}(n,1,\mR),\, \gop'=(\mathfrak{so}(n-1,\mR)\times\mR)\ltimes\mR^{n-1},
\end{eqnarray}
such that the opposite nilradicals are given by 
$\gon_-'\simeq\mR^{n-1}\subset\gon_-\simeq\mR^n$, and the one-dimensional complement 
of $\gon'_-$ in $\gon_-$ is generated by the lowest root space of $\gon_-$ (i.e., the lowest 
root space of $\gog$.) In particular, the nilradicals $\gon_-,\gon_-'$ are commutative.
The Iwasawa-Langlands decomposition of $\mathfrak{so}(n+1,1,\mR)$ is realized by the block 
decomposition
\begin{equation}
\left(
\begin{array}{ccc}
a & Y & 0 \\
X & A & -Y \\
0 & -X & -a 
\end{array}
\right),\quad a\in\mR,\, A\in \mathfrak{so}(n,\mR),\, Y\in\gon, X\in\gon_- .
\end{equation}
Let us consider the family of scalar generalized Verma $\gog$-modules 
$M^\gog_\gop(\lambda)$
induced from complex characters $\xi_\lambda: \gop\to\mC$, $\lambda\in\mC$.
The branching problem for the pair $\gog, \gog'$ applied to this class of modules
was solved in \cite{koss}, and we have
\begin{eqnarray}\label{confparbranch}
M^\gog_\gop(\lambda)|_{(\gog',\gop')}\simeq \bigoplus\limits_{j=0}^\infty M^{\gog'}_{\gop'}(\lambda-j)
\end{eqnarray}
 in the Grothendieck group $K({\fam2 O}^{\gop'})$ of the
BGG parabolic category ${\fam2 O}^{\gop'}$. 
In what follows we construct an $\mathfrak{sl}(2,\mC)$-module structure on 
$\gog'$-singular vectors (the generators of the $\gog'$-modules on the 
right hand side of \eqref{confparbranch}) in the family of $\gog$-modules 
$M^{\gog}_{\gop}(\lambda)$. 

Let $C^\alpha_l(x)$ be the $l$-th Gegenbauer polynomial in the variable $x$ with spectral 
parameter $\alpha\in\mC$, $l\in\mN_0$; we also set $C^\alpha_{-1}(x)=0$.  
See the Appendix for basic properties of Gegenbauer polynomials. 
The recurrence relations for Gegenbauer polynomials imply
\begin{eqnarray}\label{sl2gegen}
& & \big((1-x^2)\partial_x+lx\big)C^\alpha_l(x)=(l+2\alpha-1)C^\alpha_{l-1}(x),
\nonumber \\
& & \big((1-x^2)\partial_x-(l+2\alpha)x\big)C^\alpha_l(x)=-(l+1)C^\alpha_{l+1}(x).
\end{eqnarray}
\begin{lemma}  
Let $\alpha\in\mC$, $l\in\mN_0$, and let $\mathrm{Deg}$ 
be the degree operator acting on a polynomial 
of degree $l$ with the eigenvalue $l$. The linear operators 
$\{e(l),f(l),h(l))\}_{l\in\mN_0}$,
\begin{eqnarray} 
& & e(l)=(1-x^2)\partial_x-(l+2\alpha)x : \langle C^\alpha_l(x)\rangle
\,\to\, \langle C^\alpha_{l+1}(x)\rangle,
\nonumber \\
 && f(l)=(1-x^2)\partial_x+lx : \langle C^\alpha_l(x)\rangle
	\,\to\, \langle C^\alpha_{l-1}(x)\rangle ,
\nonumber \\
  &&h(l)=2(\mathrm{Deg}+\alpha) : \langle C^\alpha_l(x)\rangle
	\,\to\, \langle C^\alpha_{l}(x)\rangle
\end{eqnarray}
act on the vector 
space spanned by Gegenbauer polynomials $\{C^\alpha_l(x)\}_{l\in\mN_0}$ and furnish
it with the structure of a lowest weight $\mathfrak{sl}(2,\mC)$-module.
\end{lemma}
{\bf Proof:}
By direct computation, we have for $l\in\mN_0$:
\begin{eqnarray}
 [e,f]C^\alpha_l(x) &=& (e(l-1)f(l)-f(l+1)e(l))C^\alpha_l(x) \nonumber \\ \nonumber
&=& 2(l+\alpha)C^\alpha_l(x)=h(l)C^\alpha_l(x),
\nonumber \\
\text{and}\quad [h,e]C^\alpha_l(x) &=& (h(l+1)e(l)-e(l)h(l))C^\alpha_l(x)
\nonumber \\
&=& 2e(l)C^\alpha_l(x).
\end{eqnarray}
The same computation applies to the commutator $[h,f]$.
Notice that the collection of all $e(l)$ defines an operator $e$,
and the same for $h$ and $f$.

\hfill
$\square$

Let us now briefly review the relation of the set of singular vectors 
generating the $\gog'$-submodules on the right hand side of 
\eqref{confparbranch} to Gegenbauer polynomials, \cite{koss}. 
Denoting by $\xi_1,\dots ,\xi_{n-1},\xi_n$
the Fourier transforms of the root spaces in $\gon_-$ such that 
$\xi_1,\dots ,\xi_{n-1}$ correspond to
$\gon_-'$, $M^\gog_\gop(\lambda)\simeq \mC[\xi_1, \ldots ,\xi_{n-1}, \xi_n]$
as a vector space and the $\gog'$-singular vectors are given by homogeneous polynomials
$$
\tilde{F}_l(\xi_1,\dots ,\xi_{n-1},\xi_n)=\xi_n^l\tilde{C}^{\alpha}_l(-t^{-1}),
$$
where $l\in\mN_0$ denotes the homogeneity of the polynomial, $\alpha=-\lambda-\frac{n-1}{2}$, 
$t=\frac{1}{\xi_n^2}{\sum\limits_{j=1}^{n-1}\xi_j^2}$ and $\tilde{C}^{\alpha}_l(-t^{-1})$ is 
defined as follows: due to the fact that 
$x^{-l}C^\alpha_l(x)$ is an even rational function, we define
$x^{-l}C^\alpha_l(x)=\tilde{C}^{\alpha}_l(x^2)=\tilde{C}^{\alpha}_l(-t^{-1})$ with 
$x^2=-t^{-1}$. The space of all singular vectors is exhausted by $l\in\mN_0$.
\begin{example}
In what follows we use the normalized singular vectors $F_l(\xi_1,\dots ,\xi_{n-1},\xi_n)$, 
whose coefficient by the highest power of the quadratic invariant 
$\sum\limits_{i=1}^{n-1}\xi_i^2$ is $\lambda$-independent: 
\begin{eqnarray}
& & F_0(\xi_1,\dots ,\xi_{n-1},\xi_n)=1,
\nonumber \\
& & F_1(\xi_1,\dots ,\xi_{n-1},\xi_n)=\xi_n,
\nonumber \\
& & F_2(\xi_1,\dots ,\xi_{n-1},\xi_n)=-(2\lambda+n-3)\xi_n^2+\sum_{i=1}^{n-1}\xi_i^2,
\nonumber \\
& & F_3(\xi_1,\dots ,\xi_{n-1},\xi_n)=-(2\lambda+n-5)\xi_n^3+3\xi_n\sum_{i=1}^{n-1}\xi_i^2,
\end{eqnarray}
and so 
\begin{eqnarray}
& & \tilde{C}^\alpha_0(-t^{-1})=1,
\nonumber \\
& & \tilde{C}^\alpha_1(-t^{-1})=2\alpha,
\nonumber \\
& & \tilde{C}^\alpha_2(-t^{-1})=\alpha(t+2(1+\alpha)),
\nonumber \\
& & \tilde{C}^\alpha_3(-t^{-1})=\frac{2}{3}\alpha(\alpha+1)(3t+2(2+\alpha)),
\end{eqnarray}
where the Gegenbauer polynomials are 
\begin{eqnarray}
& & {C}^\alpha_0(x)=1,
\nonumber \\
& & {C}^\alpha_1(x)=2\alpha x,
\nonumber \\
& & {C}^\alpha_2(x)=-\alpha+2\alpha(1+\alpha)x^2,
\nonumber \\
& & {C}^\alpha_3(x)=-2\alpha(1+\alpha)x+\frac{4}{3}\alpha(1+\alpha)(2+\alpha)x^3.
\end{eqnarray}
\end{example}
The relation between $C^\alpha_l(x)$ and $\tilde{C}^\alpha_l(-t^{-1})$ implies that 
the operator identities in \eqref{sl2gegen} transform in the variable $t$ into  
\begin{eqnarray}
& & \big(-2(t+1)\partial_t+l\big)\tilde{C}^\alpha_l(-t^{-1})
=(l+2\alpha-1)\tilde{C}^\alpha_{l-1}(-t^{-1}),
\nonumber \\
& & \big(2t(t+1)\partial_t-lt-2(l+\alpha)\big)\tilde{C}^\alpha_l(-t^{-1})
=-(l+1)\tilde{C}^\alpha_{l+1}(-t^{-1}).\nonumber \\
\end{eqnarray}
The proof of the following claim is an elementary consequence of
the commutativity of the nilradical $\gon$.
\begin{lemma}
Let $\square^\xi=\sum\limits_{i=1}^{n}\partial^2_{\xi_i}$, 
$\mE_\xi=\sum\limits_{i=1}^n\xi_i\partial_{\xi_i}$,
$\alpha=-\lambda-\frac{n-1}{2}$, $l\in\mN_0$.  
Then the root space in $\gon/\gon'$ acts on the generalized Verma module 
$M^\gog_\gop(\lambda)$ by the operator 
$P(\lambda)\equiv P_n(\lambda)=i(\frac{1}{2}\xi_n\square^\xi+(\lambda-\mE_\xi)\partial_{\xi_n})$,
and descends to the map
\begin{eqnarray}
P(\lambda) :\, 
\langle F_l(\xi_1,\dots ,\xi_{n-1},\xi_n)\rangle \,
\to\, \langle F_{l-1}(\xi_1,\dots ,\xi_{n-1},\xi_n)\rangle .\nonumber 
\end{eqnarray}
In the variable $t$, $P(\lambda)$ acts by the operator
\begin{eqnarray}
\big(-2(t+1)\partial_t+l\big):\, \langle\tilde{C}^\alpha_l(-t^{-1})\rangle\,\to \, 
\langle\tilde{C}^\alpha_{l-1}(-t^{-1})\rangle.
\end{eqnarray}
\end{lemma}

Let us notice that due to the fact that the annihilator ideal of $\gog'$-singular vectors
in the algebraic Weyl algebra on $\mC^n$ is rather large, there is a plenty of algebraic 
differential operators 
inducing linear action on the vector space of $\gog'$-singular vectors (e.g., the 
same as the operator $P(\lambda)$.) In general, we can not expect to get 
$\mathfrak{sl}(2,\mC)$-actions staying entirely inside $U(\mathfrak{g})$. 
We show an example demonstrating this phenomenon that comes as close to this as
possible but still fails by considering the operator
\begin{eqnarray}
Q(\lambda):=(\sum_{i=1}^{n-1}\xi_i^2)(\frac{1}{2}\xi_n\square^\xi+(\lambda-\mE_\xi)\partial_{\xi_n})
-(\lambda-\mE_\xi+2)(n+2\lambda-2\mE_\xi+1)\xi_n,
\end{eqnarray}
which induces an action on singular vectors $F_l(\xi_1,\dots ,\xi_n)$, $l\in\mN_0$,
\begin{eqnarray}
Q(\lambda):\, \langle F_l(\xi_1,\dots ,\xi_n)\rangle\,\to\, \langle F_{l+1}(\xi_1,\dots ,\xi_n)\rangle .
\end{eqnarray}
In particular, we have 
\begin{eqnarray}
& & Q(\lambda)(1)=-(\lambda+1)(n+2\lambda-1)\xi_n, \nonumber \\
& & Q(\lambda)(\xi_n)=\lambda\big((\sum_{i=1}^{n-1}\xi_i^2)-(n+2\lambda-3)\xi_n^2\big), \nonumber \\
& & Q(\lambda)((\sum_{i=1}^{n-1}\xi_i^2)-(n+2\lambda-3)\xi_n^2)
=-(\lambda-1)(n+2\lambda-3)\big(3\xi_n(\sum_{i=1}^{n-1}\xi_i^2)-(n+2\lambda-5)\xi_n^3\big), \nonumber \\
& & Q(\lambda)(3\xi_n(\sum_{i=1}^{n-1}\xi_i^2)-(n+2\lambda-5)\xi_n^3)=(\lambda-2)
\big(3(\sum_{i=1}^{n-1}\xi_i^2)^2 \nonumber \\
& & -6(n+2\lambda-5)\xi_n^2(\sum_{i=1}^{n-1}\xi_i^2)
+(n+2\lambda-7)(n+2\lambda-5)\xi_n^4\big). \nonumber \\
\end{eqnarray}
Notice that the operator $Q(\lambda)$ is an 
element of the universal enveloping algebra ${U}(\gog)$. 
A disadvantage of $Q(\lambda)$ is that the pair of operators 
$P(\lambda), Q(\lambda)$ together with the homogeneity operator do 
not close in an $\mathfrak{sl}(2,\mC)$-algebra realized in
${ U}(\gog)$. A straightforward but tedious computation shows
\begin{eqnarray}
& & [\frac{1}{2}\xi_n\square^\xi+(\lambda-\mE_\xi)\partial_{\xi_n}, Q(\lambda)]=
 -\frac{1}{2}(4\lambda+10+2\mE_\xi)\xi_n^2\square^\xi \nonumber \\
& & +[(2\mE_\xi+n-3)(\lambda+1-\mE_\xi) +(\lambda-\mE_\xi+2)(n+2\lambda+1-2\mE_\xi) \nonumber \\
& & -(\lambda-\mE_\xi)(n+4\lambda+3-4\mE_\xi)-(n+4\lambda+7+4\mE_\xi)]\xi_n\partial_{\xi_n} \nonumber \\
& & -((\sum_{i=1}^{n-1}\xi_i^2)+\xi_n^2)(\xi_n\partial_{\xi_n}+1)\square^\xi 
-2(\lambda+1-\mE_\xi)((\sum_{i=1}^{n-1}\xi_i^2)+\xi_n^2)\partial_{\xi_n}^2 \nonumber \\
& & -(\lambda-\mE_\xi)[(\lambda-\mE_\xi+2)(n+2\lambda+1-2\mE_\xi)-(n+4\lambda+3-4\mE_\xi)],
\nonumber \\
\end{eqnarray}
where the operator on the right hand side of the last equality
does not act in the
algebraic Weyl algebra on $\mC^n$ as a multiple of identity on both
$P(\lambda)$ and $Q(\lambda)$.

In what follows we construct two operators in the variables 
$\xi_1,\dots,\xi_n$, which are not the elements of  
${U}(\gog)$, but they fulfill 
$\mathfrak{sl}(2,\mC)$-commutation relations 
when their action is restricted to the set of $\gog'$-singular vectors
in $M^\gog_\gop(\lambda)$.    
\begin{theorem}\label{sl2pair}
Let $l\in\mN_0$. Then the collection of operators $\{e_\xi(l), f_\xi(l), h_\xi(l)\}_{l\in\mN_0}$
in the variables $\xi_1,\dots ,\xi_n$,
\begin{eqnarray}
& & e_\xi(l):=-(\sum_{i=1}^n\xi_i^2){\partial_{\xi_n}}-(l+2\alpha)\xi_n:\,
\langle F_l(\xi_1,\dots, \xi_n)\rangle\,\to\, \langle F_{l+1}(\xi_1,\dots, \xi_n)\rangle,
\nonumber \\
& & f_\xi(l):=\frac{1}{\xi_n}\big(\frac{(\sum_{i=1}^n\xi_i^2)}{(\sum_{i=1}^{n-1}\xi_i^2)}
(\xi_n{\partial_{\xi_n}}-l) +l\big):\,
\langle F_l(\xi_1,\dots, \xi_n)\rangle\,\to\, \langle F_{l-1}(\xi_1,\dots, \xi_n)\rangle,
\nonumber \\ \label{fxil}
& & h_\xi(l):=2(l+\alpha)\mathrm{Id} :\,
\langle F_l(\xi_1,\dots, \xi_n)\rangle\,\to\, \langle F_{l}(\xi_1,\dots, \xi_n)\rangle,
\nonumber \\
\end{eqnarray} 
fulfill the $\mathfrak{sl}(2,\mC)$-commutation relations, and 
$F_l(\xi_1,\dots ,\xi_n), l\in\mN_0$, are the one-dimensional 
weight spaces of a highest weight ${\mathfrak sl}(2,\mC)$-Verma module.
\end{theorem}
{\bf Proof:}
Let $f=f(\xi_1,\dots, \xi_n)$, $l\in\mN_0$, 
and $t=\frac{1}{\xi_n^2}{\sum\limits_{i=1}^{n-1}\xi_i^2}$. We have 
$$
{\partial_{\xi_n}}(\xi_n^{-l}f)=-l\xi_n^{-l-1}f+\xi_n^{-l}{\partial_{\xi_n}}f,
$$
and 
$$
\partial_t=-\frac{\xi_n^3}{2(\sum\limits_{i=1}^n\xi_i^2)}{\partial_{\xi_n}}=
\frac{\xi_n^2}{2\xi_i}{\partial_{\xi_i}}
$$
for all $i=1,\dots ,n-1$. By direct substitution for $t$, the first operator equals to 
\begin{eqnarray}
\xi_n^{l+1}\big(
2\frac{(\sum\limits_{i=1}^{n-1}\xi_i^2)}{\xi_n^2}(\frac{(\sum\limits_{i=1}^{n-1}\xi_i^2)}{\xi_n^2}+1)
(\frac{-\xi_n^3}{2(\sum\limits_{i=1}^{n-1}\xi_i^2)}{\partial_{\xi_n}})
-l\frac{(\sum\limits_{i=1}^{n-1}\xi_i^2)}{\xi_n^{2}}-2(\alpha+l)
\big)\xi_n^{-l},\nonumber \\
\end{eqnarray}
and standard manipulations give the required result $e_\xi(l)$. Analogously, the second operator 
is equal to 
\begin{eqnarray}
\xi_n^{l-1}\big(
-2(\frac{(\sum\limits_{i=1}^{n-1}\xi_i^2)}{\xi_n^2}+1)
(\frac{-\xi_n^3}{2(\sum\limits_{i=1}^{n-1}\xi_i^2)}{\partial_{\xi_n}})
+l
\big)\xi_n^{-l},
\end{eqnarray}
and this gives $f_\xi(l)$.

As for the $\mathfrak{sl}(2,\mC)$-commutation relations, 
it is again straightforward to check that
\begin{eqnarray}
& & f_\xi(l+1)e_\xi(l)-e_\xi(l-1)f_\xi(l)= \nonumber \\
& &  [\frac{1}{\xi_n}\big(\frac{(\sum\limits_{i=1}^n\xi_i^2)}{(\sum\limits_{i=1}^{n-1}\xi_i^2)}
    (\xi_n{\partial_{\xi_n}}-(l+1)) +(l+1)\big)] \circ
    [-(\sum\limits_{i=1}^n\xi_i^2){\partial_{\xi_n}}-(l+2\alpha)\xi_n]
    \nonumber \\
& & -[-(\sum\limits_{i=1}^n\xi_i^2){\partial_{\xi_n}}-(l-1+2\alpha)\xi_n]\circ
    [\frac{1}{\xi_n}\big(\frac{(\sum\limits_{i=1}^n\xi_i^2)}{(\sum\limits_{i=1}^{n-1}\xi_i^2)}
     (\xi_n{\partial_{\xi_n}}-l) +l\big)] \nonumber \\
\end{eqnarray}
is equal to $-h_\xi(l)$. The two remaining relations involving  $h_\xi(l)$ are the consequence of the 
explicit formula for $h_\xi(l)$ and the homogeneity of $e_\xi(l)$, $f_\xi(l)$.

Finally, the structure of $\mathfrak{sl}(2,\mC)$-Verma module 
is a consequence of $\mathfrak{sl}(2,\mC)$-commutation 
relations and the action of $e_\xi(l), f_\xi(l)$ and $h_\xi(l)$ on the weight vectors 
$F_l(\xi_1,\dots ,\xi_n)$.

\hfill
$\square$

\begin{remark}
The normalization of singular vectors $F_l(\xi_1,\dots,\xi_n)$, $l\in\mN_0$,
is chosen in such a way that the action of particular basis elements 
$e_\xi(l), f_\xi(l)$ produces the constants:

\vspace{0.5cm}
\xymatrix{
& & & & F_{2l} \ar@/_/[d]_{-(2\alpha+2l)} & & \\
& & & & F_{2l+1} \ar@/_/[d]_{-1} \ar@/_/[u]_{2l+1} & &\quad ,  \\
& & & & F_{2l+2} \ar@/_/[u]_{(2l+2)(2\alpha+2l+1)} & &   \\
& & & & & & \\
}

and their commutator gives the weight $2\alpha+4l-2$. 
\end{remark}
The operator 
$\frac{1}{2}\xi_n\square^\xi+(\lambda-\mE_\xi)\partial_{\xi_n}$
does not have $\mathfrak{sl}(2,\mC)$-commutation 
relation with $e_\xi(l):=-(\sum_{i=1}^n\xi_i^2){\partial_{\xi_n}}-(l+2\alpha)\xi_n$.
In fact, in the algebraic Weyl algebra we have 
\begin{eqnarray}
& & [-(\sum_{i=1}^n\xi_i^2){\partial_{\xi_n}}-(\mE_\xi-1+2\alpha)\xi_n,
\frac{1}{2}\xi_n\square^\xi+(\lambda-\mE_\xi)\partial_{\xi_n}]= \nonumber \\
& & -\frac{1}{2}(\sum_{i=1}^{n-1}\xi_i^2)\square^\xi
-(\sum_{i=1}^{n-1}\xi_i^2){\partial_{\xi_n}}{\partial_{\xi_n}}
+\frac{1}{2}\xi_n^2\square^\xi  \nonumber \\
& & +(n+\lambda+\mE_\xi)\xi_n{\partial_{\xi_n}}+(\mE_\xi+2\alpha)(\lambda-\mE_\xi).
\end{eqnarray}

Notice that the operator $f_\xi(l)$ introduced in \eqref{fxil} does not 
belong to the universal enveloping algebra $U(\gog)$, but rather to its
localization with respect to the subalgebra of invariants in $U(\gon_-)$
with respect to the simple part of the Levi factor $\gol'$.

Let us finally examine the action of the $\mathfrak{sl}(2,\mC)$-Casimir 
operator on $\gog'$-singular vectors. Recall that for the Lie 
algebra $\mathfrak{sl}(2,\mC)$ generated by elements $e,f,h$ with 
commutation relations $[e,f]=h, [h,e]=2e$
and $[h,f]=-2f$, the Casimir operator is 
$\mathrm{Cas}=ef+fe+\frac{1}{2}h^2$.
\begin{theorem}
\begin{enumerate}
\item
The Casimir operator of the Lie algebra $\mathfrak{sl}(2,\mC)$ realized by 
$e_\xi(l),f_\xi(l),h_\xi(l)$, $l\in\mN_0$, is 
\begin{eqnarray}\label{casimirfinal}
& & \mathrm{Cas}:= f_\xi(l+1)e_\xi(l)+e_\xi(l-1)f_\xi(l)+\frac{1}{2}h^2_\xi(l)=
\nonumber \\
& & -2\frac{(\sum\limits_{i=1}^{n}\xi_i^2)^2}{(\sum\limits_{i=1}^{n-1}\xi_i^2)}{\partial^2_{\xi_n}}
-2(2\alpha+1)\frac{(\sum\limits_{i=1}^{n}\xi_i^2)}{(\sum\limits_{i=1}^{n-1}\xi_i^2)}\xi_n{\partial_{\xi_n}} 
\nonumber \\
& & -2\frac{\alpha(\sum\limits_{i=1}^{n-1}\xi_i^2)-l(l+2\alpha)\xi_n^2}{(\sum\limits_{i=1}^{n-1}\xi_i^2)}
+2(\mE_\xi+\alpha)^2
\end{eqnarray}
when acting on $F_l(\xi_1,\dots ,\xi_n)$.
\item
The Casimir operator $\mathrm{Cas}$ acts by $2\alpha(\alpha-1)$-multiple of identity on
the $\mathfrak{sl}(2,\mC)$-Verma module with weight spaces 
$\{\langle F_l(\xi_1,\dots ,\xi_n)\rangle \}_{l\in\mN_0}$.  
\end{enumerate}
\end{theorem}
{\bf Proof:}
We have 
\begin{eqnarray}
& & \mathrm{Cas} =f_\xi(l+1)e_\xi(l)+e_\xi(l-1)f_\xi(l)+\frac{1}{2}h^2_\xi(l)= \nonumber \\
& &  [\frac{1}{\xi_n}\big(\frac{(\sum\limits_{i=1}^n\xi_i^2)}{(\sum\limits_{i=1}^{n-1}\xi_i^2)}
    (\xi_n{\partial_{\xi_n}}-(l+1)) +(l+1)\big)] \circ
    [-(\sum\limits_{i=1}^n\xi_i^2){\partial_{\xi_n}}-(l+2\alpha)\xi_n]
    \nonumber \\
& & +[-(\sum\limits_{i=1}^n\xi_i^2){\partial_{\xi_n}}-(l-1+2\alpha)\xi_n]\circ
    [\frac{1}{\xi_n}\big(\frac{(\sum\limits_{i=1}^n\xi_i^2)}{(\sum\limits_{i=1}^{n-1}\xi_i^2)}
     (\xi_n{\partial_{\xi_n}}-l) +l\big)] \nonumber \\
& & +\frac{1}{2}(-2(\mE_\xi+\alpha))^2.     
\end{eqnarray}
Expanding the compositions and recollecting all terms according to the power 
of $\partial_{\xi_n}$, we arrive at \eqref{casimirfinal}. 
As for the proof of the second claim, it follows from 
$$
\mathrm{Cas}(1)=(-2\alpha+2\alpha^2)1=2\alpha(\alpha-1)1.
$$
The proof is complete.

\hfill
$\square$


\subsection{The pair of Lie algebras 
($\mathfrak{sl}(2,\mR)\times \mathfrak{sl}(2,\mR), \mathrm{diag}(\mathfrak{sl}(2,\mR))$)
and the Jacobi polynomials}
\label{subsection2.2}
Let us consider the pair of compatible Lie algebras and their Borel subalgebras, 
\begin{eqnarray}
& & \gog=\mathfrak{sl}(2,\mR)\times \mathfrak{sl}(2,\mR),\, 
\gog\supset\gob=(\mR\times\mR)\ltimes(\mR\times\mR),
\nonumber \\
& & \gog'=\mathrm{diag}(\mathfrak{sl}(2,\mR)),\, \gog'\supset\gob'=
\mathrm{diag}\big((\mR\times\mR)\ltimes(\mR\times\mR)\big),
\end{eqnarray}
such that $\gon_-'\subset\gon_-$ ($\mathrm{dim}_\mR(\gon_-)=2, 
\mathrm{dim}_\mR(\gon_-')=1$) and its one 
dimensional complement is denoted by $F$. 
In particular, the nilradicals $\gon_-,\gon_-'$ are commutative and $X$ denotes 
the generator of $\gon_-'$.

Let us consider the family of scalar Verma $\gog$-modules 
$M^\gog_\gop(\lambda, \mu)$, induced from complex 
characters $\xi_{\lambda, \mu}: \gob\to\mC$ with $\lambda, \mu\in\mC$.
The branching of scalar Verma modules for 
the pair $(\gog, \gob)$ and $(\gog', \gob')$  
 is given by
\begin{eqnarray}\label{sl2diagvm}
M^\gog_\gop(\lambda, \mu)|_{(\gog',\gop')}\simeq 
\bigoplus_{j=0}^\infty M^{\gog'}_{\gop'}(\lambda+\mu-2j)
\end{eqnarray}
in the Grothendieck group $K({\fam2 O}(\gog))$ 
of the BGG category ${\fam2 O}(\gog)$.

In \cite{koss2}, the generators of $\gog'$-submodules generating the 
summands on the right 
hand side of \eqref{sl2diagvm} are determined and the 
non-trivial composition structure related to factorization 
properties of Jacobi polynomials is discussed in the non-generic case 
$\lambda+\mu\in\mN_0$. Here we restrict to the case $\lambda,\mu\notin\mN_0$.
For $\nu \in {\mathbb{N}}_0$, 
we define (see \cite[Chapter 3]{hum}) 
a ${\mathfrak {g}'}$-module
 $P_{{\mathfrak {b}'}}^{\mathfrak {g}'}(\nu)$
 as the non-split extension
\begin{equation}
\label{eqn:proj}
  0 
  \to 
  M_{{\mathfrak {b}}'}^{\mathfrak {g}'}(\nu)
  \to 
  P_{{\mathfrak {b}}'}^{\mathfrak {g}'}(\nu)
  \to 
  M_{{\mathfrak {b}'}}^{\mathfrak {g}'}(-\nu-2)
  \to 0.
\end{equation}
We introduce an involution
\[
\iota: \mathbb{C} \to \mathbb{C},
\ 
\nu \mapsto -\nu-2.
\]
and for $N\in\mathbb{N}_0$ set
\begin{alignat*}{3}
&\Lambda &&\equiv \Lambda(N)
         &&:= \{N-2l: l\in\mathbb{N}_0\},
\\
&\Lambda_s &&\equiv \Lambda_s(N)
           &&:= \{N-2l: l\in\mathbb{N}_0, \  2l\le \lambda+\mu\},
\\
&\Lambda_r &&\equiv \Lambda_r(N)
           &&:= \Lambda\setminus(\Lambda_s\cup\iota(\Lambda_s))
\\
&  &&      &&\ = \begin{cases}
                 \{-1\}\cup\{-N,-N-2,-N-4,\dotsc\}
                 &\text{($N$: even)},\\
                 \{-N,-N-2,-N-4,\dotsc\}
                 &\text{($N$: odd)}.
              \end{cases}
\end{alignat*}
\begin{theorem}\cite{koss2}
\label{thm:7.5}
Suppose that $\lambda +\mu\in\mN_0$ with $\lambda,\mu\notin\mN_0$, i.e.,
the scalar Verma modules $M^{\mathfrak{sl}(2,\mR)}_{\gob'}(\lambda)$ resp. 
$M^{\mathfrak{sl}(2,\mR)}_{\gob'}(\mu)$
are irreducible $\gog'$-modules. 
Then the tensor product of two scalar Verma modules 
$
M^{\mathfrak{sl}(2)}_{\gob'}(\lambda)\otimes M^{\mathfrak{sl}(2)}_{\gob'}(\mu)$
decomposes as 
$\mathrm{diag}(\mathfrak{sl}(2,{\mathbb R}))\simeq \mathfrak{sl}(2,{\mathbb R})$-module 
\begin{equation}
 \bigoplus_{\nu\in \Lambda_r(\lambda + \mu)}M^{\mathfrak{sl}(2)}_{\gob'}(\nu)
\,\,\oplus\,\,
\bigoplus_{\nu \in \Lambda_s(\lambda +\mu)}P^{\mathfrak{sl}(2)}_{\gob'}(\nu).
\end{equation}
Here $P^{\mathfrak{sl}(2)}_{\gob'}(\nu)$ 
are the projective objects  defined in \eqref{eqn:proj}.  
\end{theorem}
It is straightforward to see that the generator $X$ 
of $\gon_+'$ acts 
in the non-compact model of the representation 
${\lC^\infty}(\gon_-, {\mathbb C}_{\lambda ,\mu})$, 
induced from the character $({\lambda ,\mu})$ of $\gob$
on the $1$-dimensional vector space ${\mathbb C}_{\lambda ,\mu}\simeq \mC$, 
by the first order differential operator
\begin{eqnarray}
d\pi(X)=\lambda x+x^2\partial_x+ \mu y+y^2\partial_y.
\end{eqnarray}
The action on the scalar Verma module, 
induced from the dual representation to 
${\mathbb C}_{\lambda ,\mu}$ and realized in the
Fourier dual picture, results into the 
second order differential operator 
\begin{eqnarray}\label{rancohop1}
d\tilde{\pi}(X)=i(-\lambda \partial_\xi+\xi\partial_\xi^2-\mu \partial_\eta+\eta\partial_\eta^2)
\end{eqnarray}
acting on polynomial algebra $\mC[\xi,\eta]$.
As a $\gog$-module, $\mC[\xi,\eta]$ can be identified
with the Verma module induced from the dual representation to 
${\mathbb C}_{\lambda ,\mu}$.
Analogously, the generator $F$ of $\gon_+/\gon_+'$ acts 
on ${\lC^\infty}(\gon_-, {\mathbb C}_{\lambda ,\mu})$ and $\mC[\xi,\eta]$
by differential operators
\begin{eqnarray}
d\pi(F)=\lambda x+x^2\partial_x- \mu y-y^2\partial_y
\end{eqnarray}
and 
\begin{eqnarray}\label{rancohop2}
d\tilde{\pi}(F)=i(-\lambda \partial_\xi+\xi\partial_\xi^2+\mu \partial_\eta-\eta\partial_\eta^2),
\end{eqnarray}
respectively.
The Levi factor $\gol\subset\gob$ contains the Euler homogeneity operator and the operator  
$d\tilde{\pi}(X)$ preserves the space of homogeneous polynomials. We define
$t=\frac{\xi}{\eta}$ and write a homogeneous polynomial as $\eta^l Q(t)$ for some polynomial $Q=Q(t)$
 of
degree $l$. We easily compute
\begin{eqnarray}
 \partial_{\xi}(\eta^l Q(t))&=&\eta^{l-1} Q^\prime(t),
\nonumber 
\\
 \partial_{\eta}(\eta^l Q(t))
 &=&
 \eta^{l-1}l Q-\xi \eta^{l-2}Q^\prime(t),
\nonumber 
\\
 \partial^2_{\xi}(\eta^l Q(t))
 &=&
 \eta^{l-2}Q^{\prime\prime}(t),
\nonumber \\
 \partial^2_{\eta}(\eta^l Q(t))
 &=&\eta^{l-2}l(l-1)Q(t)-2(l-1)\xi\eta^{l-3} Q^\prime(t)
\nonumber \\
& & +\xi^2 \eta^{l-4} Q^{\prime\prime}(t).  
\label{eqn:7.7}
\end{eqnarray}
The substitution into \eqref{rancohop1}
 resp. \eqref{rancohop2} yields the differential equation
consisting of polynomials 
 $\eta^l Q(\frac \xi \eta)$,
 where $Q(t)$ is a polynomial solution to  
\begin{equation}
\label{hypergeomeq1}
[t(t+1) \partial_t^2
 + (t(\mu-2(l-1))-\lambda)\partial_t
+l(l-1-\mu)]Q(t)=0,
\end{equation}
resp. the differential equation representing the action of $F\in \gon_+/\gon'_+$: 
\begin{equation}
\label{hypergeomeq2}
-t(t-1) \partial_t^2
 + (t(2l-\mu-2)-\lambda)\partial_t
+l(\mu-l+1).
\end{equation} 
The ordinary second order hypergeometric differential equation (\ref{hypergeomeq1})
is the Jacobi differential equation, and its polynomial solutions are the Jacobi
polynomials:
 \begin{theorem}\cite{koss2}
\label{thm:SLSol}
Let $P^{\alpha,\beta}_l(x)$ denote the degree $l$ polynomial solution
of the Jacobi hypergeometric equation \eqref{eqn:9.3}, see \eqref{jacpol}.
Let us define homogeneous polynomials $\tilde{P}_{l}^{-\lambda-1, \mu + \lambda -2l+1}(\xi,\eta)$ by
\begin{eqnarray}
\tilde{P}_{l}^{-\lambda-1, \mu + \lambda -2l+1}(\xi,\eta):=
\eta^l P_{l}^{-\lambda-1, \mu + \lambda -2l+1}\left(\frac{2 \xi}{\eta}+1\right)
\end{eqnarray}
for all $l\in\mN_0$.
Then for any $\lambda, \mu \in {\mathbb{C}}$, 
\begin{eqnarray}
  \bigoplus_{l=0}^{\infty}\langle\tilde{P}_{l}^{-\lambda-1, \mu + \lambda -2l+1}\rangle 
\end{eqnarray}
is the complete set of polynomial solutions of (\ref{rancohop1}) representing the 
singular vectors in the Fourier dual picture.
\end{theorem}
By abuse of notation, we denote by $d\tilde{\pi}(X)$, $d\tilde{\pi}(F)$ the two 
mutually commuting operators
acting on degree $l$-polynomials in the non-homogeneous variable $t$:
\begin{eqnarray}
& & d\tilde{\pi}(X):=t(t+1) \partial_t^2
 + (t(\mu-2(l-1))-\lambda)\partial_t
+l(l-1-\mu),
\nonumber \\
& & d\tilde{\pi}(F):=-t(t-1) \partial_t^2
 + (t(2l-\mu-2)-\lambda)\partial_t
+l(\mu-l+1).
\end{eqnarray}
\begin{theorem}
Let $P_l^{\alpha,\beta}(x)$ be the Jacobi polynomial of degree $l\in\mN_0$ and  
$$
d\tilde{\pi}(X)(P_l^{-\lambda-1, \mu + \lambda -2l+1}(2t+1))=0
$$ 
for $x=2t+1$. Then 
\begin{eqnarray}
& & d\tilde{\pi}(F)(P_l^{-\lambda-1, \mu + \lambda -2l+1}(2t+1))=
\nonumber \\
& & 2(l-1-\lambda)(\mu-l+1) P_{l-1}^{-\lambda-1, \mu + \lambda -2l+3}(2t+1),
\end{eqnarray}  
i.e., $d\tilde{\pi}(F)$ maps the homogeneity 
$l$ polynomial solution of the Jacobi differential equation $d\tilde{\pi}(X)$
to (a multiple depending on $\alpha,\beta$ of)
the homogeneity $(l-1)$ polynomial solution of the Jacobi differential equation.
\end{theorem} 
{\bf Proof:}
We first observe that $\sum\limits_{i=0}^la^l_it^i$ is the degree $l$ Jacobi polynomial 
in the variable $t$ provided
the recursion relations 
\begin{eqnarray}\label{llevel}
[i(i-2l+\mu+1)+l(l-\mu-1)]a^l_i+(i-\lambda)(i+1)a^l_{i+1}=0
\end{eqnarray}
are satisfied for all $i=0,\dots ,l$. 
The operator $d\tilde{\pi}(F)$ maps the degree $l$ Jacobi polynomial 
to the space of polynomials of degree $(l-1)$. The 
reason is that the coefficient of the monomial $t^l$ in 
$d\tilde{\pi}(F)(P_l^{-\lambda-1, \mu + \lambda -2l+1}(2t+1))$  is equal to
$$
-l(l-1)+l(2l-\mu-2)+l(\mu-l+1)=0.
$$ 
In particular, we have 
$$
d\tilde{\pi}(F)(P_l^{-\lambda-1, \mu + \lambda -2l+1}(2t+1))=
\sum\limits_{i=0}^l[2i(-i+2l-\mu-1)+2l(\mu-l+1)]a^l_it^i
$$
for $i=0,\dots ,l-1$, and it remains to prove that this 
polynomial is, up to a multiple,
the degree $(l-1)$ Jacobi polynomial. Assuming the recursion relation 
(\ref{llevel}) holds, we prove that 
$[2i(-i+2l-\mu-1)+2l(\mu-l+1)]a^l_i$ are the coefficients of the degree 
$l-1$ Jacobi polynomial (see again \eqref{llevel}):
\begin{eqnarray} 
& & [i(i-2l+\mu+3)+(l-1)(l-\mu-2)][2i(-i+2l-\mu-1)+
\nonumber \\
& & 2l(\mu-l+1)]a^l_i= -(i-\lambda)(i+1)[2(i+1)(-i+2l-\mu-2)+
\nonumber \\
& & 2(\mu-l+1)]a^l_{i+1}.
\end{eqnarray}
However, the last equality is equivalent to 
$$
[i(i-2l+\mu+3)+(l-1)(l-\mu-2)]=-[(i+1)(-i+2l-\mu-2)+l(\mu-l+1)],
$$
which is easy to verify and the claim follows. 

It remains to compute the explicit polynomial in $\lambda, \mu$ as 
a coefficient of the proportionality. By definition,
\begin{align}
& a^l_l=\frac{1}{l!}\mu(\mu-1)\dots (\mu-l+2)(\mu-l+1),\nonumber \\
& a^{l-1}_{l-1}=\frac{1}{(l-1)!}\mu(\mu-1)\dots (\mu-l+2), \nonumber
\end{align}
and so we get for $i=l-1$
$$
d\tilde{\pi}(F)(\sum_{i=0}^la^l_it^i)=2\mu a^{l-1}_lt^{l-1}+\cdots .
$$
Because $-\mu a_l^{l-1}=-l(l-1-\lambda)a_l^l$, a direct comparison yields 
the required form $2(l-1-\lambda)(\mu-l+1)$ of the coefficient of proportionality.
The proof is complete.
\hfill
$\square$

\begin{example}
Let us present a first few low degree polynomials: 
\begin{eqnarray}
& & P_0^{-\lambda-1, \mu + \lambda +1}(2t+1)=1,
\nonumber \\
& & P_1^{-\lambda-1, \mu + \lambda -1}(2t+1)=\mu t-\lambda,
\nonumber \\
& & P_2^{-\lambda-1, \mu + \lambda -3}(2t+1)=-\frac{1}{2}[\mu(1-\mu)t^2+2(1-\mu)(1-\lambda)t+\lambda(1-\lambda)],
\nonumber \\
& & P_3^{-\lambda-1, \mu + \lambda -5}(2t+1)=\frac{1}{6}[\mu(1-\mu)(2-\mu)t^3+3(2-\mu)(\mu-1)(2-\lambda)t^2
\nonumber \\
& & +3(2-\mu)(1-\lambda)(2-\lambda)t+\lambda(2-\lambda)(1-\lambda)].
\nonumber
\end{eqnarray}
It is straightforward to check
\begin{eqnarray}
& &
d\tilde{\pi}(F)(P_1^{-\lambda-1, \mu + \lambda -1}(2t+1))=-2\lambda\mu P_0^{-\lambda-1, \mu + \lambda +1}(2t+1),
\nonumber \\
& & 
d\tilde{\pi}(F)(P_2^{-\lambda-1, \mu + \lambda -3}(2t+1))=-2(\lambda-1)(\mu-1) P_1^{-\lambda-1, \mu + \lambda -1}(2t+1),
\nonumber \\
& & d\tilde{\pi}(F)(P_3^{-\lambda-1, \mu + \lambda -5}(2t+1))=-2(\lambda-2)(\mu-2) P_2^{-\lambda-1, \mu + \lambda -3}(2t+1).
\nonumber 
\end{eqnarray}
\end{example}
We remark that because of the assumption of irreducibility, $\lambda,\mu\notin\mN_0$,
the coefficients of proportionality $(2(l-1-\lambda)(\mu-l+1))$ are non-zero. In this 
example we do not attempt to construct the complete ${\mathfrak sl}(2,\mC)$-structure
on the space of ${\gog}'$-singular vectors. 


\section{Relative Lie and Dirac cohomology and $\gog'$-singular vectors}

The Lie algebra (co)homology or the Dirac cohomology associated to a 
Lie algebra and its modules are among important algebraic invariants with 
applications in representation theory, see \cite{koslac}, \cite{che}, \cite{hp}.
 
In fact, the two examples in Section $2$, Section $3$ are motivated by the following
general problem. Let us consider the short exact sequence of pairs 
of Lie algebras and their parabolic subalgebras:
\begin{eqnarray}
0\to (\gog',\gop')\to (\gog,\gop)\to (\gog,\gop)/(\gog',\gop')\to 0
\end{eqnarray}
and a $\gog$-module $V$. In our applications, $\gog$ and $\gog'$ $(\gog'\subset\gog)$ 
are simple Lie algebras, $\gop\subset\gog$ and $\gop'\subset\gog'$ their 
parabolic subalgebras  ($\gop$ is $\gog'$-compatible), and the vector complements of $\gop$ and $\gop'$ in $\gog$ 
and $\gog'$ are the Lie algebras of the opposite nilradicals
$\gon_-$ and $\gon_-'$.

Then the key question is an intrinsic definition of the relative Lie algebra 
(co)differential or relative Dirac operator associated to  
compatible couples of Lie algebras given by simple Lie algebra and its parabolic
subalgebra, and their role in the compatibility of the branching problem 
applied to $\gog'\subset\gog$ and the parabolic BGG category ${\fam2 O}^\gop(\gog)$. 

Namely, for $\gon_-=\gon_-'\oplus(\gon_-/\gon_-')$ with
$\gon_-'$ the ideal in $\gon$, we would like to define 
the relative Dirac operator such that the underlying 
relative Dirac cohomology functor $H_{D,\mathrm{rel}}(\gon_-,\gon_-'; -)$
abuts in a spectral sequence to the $(\gog,\gop)$-Dirac 
cohomology of $M^\gog_\gop(\lambda)$:
\begin{eqnarray}
H_D(\gog,\gop; M^\gog_\gop(\lambda))\Longrightarrow
H_{D,\mathrm{rel}}(\gon_-,\gon_-', H_D(\gog',\gop'; M^\gog_\gop(\lambda)|_{\gog'}) .
\end{eqnarray}
Here $H_D(\gog,\gop; M^\gog_\gop(\lambda))$ was determined in \cite{hs} for
irreducible generalized Verma modules, while in the presence of non-trivial 
composition series the relevant higher Dirac cohomology was constructed in 
\cite{ps}. 

We do not have an answer to the previous question, and the examples in
Section \ref{subsection2.1}, Section \ref{subsection2.2} clearly demonstrate the difficulties.
To be more explicit, we shall stick to the case discussed in Section \ref{subsection2.1} 
and assume that $\lambda\in\mC$
is generic so that $M^\gog_\gop(\lambda)$ as well as $M^{\gog'}_{\gop'}(\lambda-j)$
are irreducible highest weight modules for all $j\in\mN_0$. 
Then the $(\gog,\gop)$-Dirac cohomology 
of the left hand side \eqref{confparbranch} equals to
$$ 
H_D(\gog,\gop,M^\gog_\gop(\lambda))\simeq \mC_\lambda\otimes\mC_{\rho(\gon_-)},
$$ 
while the $(\gog', \gop')$-Dirac cohomology of the right hand side 
of \eqref{confparbranch} is 
$$
H_D(\gog',\gop',\bigoplus\limits_{j=0}^\infty M^{\gog'}_{\gop'}(\lambda-j))\simeq 
\bigoplus\limits_{j=0}^\infty\mC_{\lambda-j}\otimes\mC_{\rho(\gon_-')}. 
$$  
Here we used the notation $\rho(\gon_-)$ and $\rho(\gon_-')$ for the half sum 
of roots of root spaces in $\gon_-$ and $\gon_-'$, respectively, and 
$\mC_\lambda$ denotes 
the one dimensional inducing representation of $\xi_\lambda$.  
Then the relative Dirac operator $D_{\mathrm{rel}}$
\begin{eqnarray}\label{reldirex}
D_{\mathrm{rel}}(\gon_-,\gon_-'):=e_\xi\otimes f_\xi+f_\xi\otimes e_\xi,
\end{eqnarray}
based on Lemma \ref{sl2pair}
and the highest weight $\mathfrak{sl}(2,\mC)$-module with the action
of $\{e_\xi(l), f_\xi(l), h_\xi(l)\}_{l\in\mN_0}$,
computes the expected result. 

However, this approach clearly fails for
non-generic values of $\lambda$ because the ${\mathfrak sl}(2,\mC)$-module 
is no longer irreducible. Moreover, the construction 
\eqref{reldirex} does not intrinsically proceed in $U(\gog)$.

\section{Appendix: Jacobi and Gegenbauer polynomials}

In the present section we summarize for the reader's convenience 
a few basic conventions and properties related to the Jacobi 
and Gegenbauer polynomials.

We use the notation $\Gamma(z)$ for the Gamma function, $z\in\mC$, and the 
analytical continuation of 
the binomial coefficient 
is given by
$$
{z\choose l} := \frac{\Gamma(z+1)}{\Gamma(l+1)\Gamma(z-l+1)},\quad
{z \choose l}=0
\,\,\text{if}\,\, l-z \in {\mathbb{N}}
 \,\,\text{and}\,\, z \not\in -{\mathbb{N}}.
$$

The Jacobi polynomials $P_{l}^{(\alpha,\beta)}(z)$
of degree $l\in\mN_0$ with two spectral parameters $\alpha, \beta\in\mC$
are defined as special values of the hypergeometric function
\begin{align}
P_l^{(\alpha,\beta)}(z)
=&{l+\alpha\choose l}
\,_2F_1\left(-l,1+\alpha+\beta+l;\alpha+1;\frac{1-z}{2}\right)
\nonumber \\
=&
\frac{\Gamma (\alpha+l+1)}{l!\Gamma (\alpha+\beta+l+1)}
\sum_{m=0}^n {l\choose m}
\frac{\Gamma (\alpha + \beta + l + m + 1)}{\Gamma (\alpha + m + 1)}
\left(\frac{z-1}{2}\right)^m
\nonumber \\
=&
\sum_{j=0}^l
{l+\alpha\choose j}{l+\beta \choose l-j}
\left(\frac{z-1}{2}\right)^{l-j} \left(\frac{z+1}{2}\right)^{j}, 
\label{jacpol}
\end{align}
normalized by
$$
P_l^{(\alpha, \beta)} (1) = {l+\alpha\choose l}=\frac{(\alpha+1)_l}{l!}.
$$
Here $(\alpha+1)_l=(\alpha+1)(\alpha+2)\cdots (\alpha+l)$ 
denotes the Pochhammer symbol for the partially rising factorial.
The Jacobi polynomials satisfy the orthogonality relations
\begin{eqnarray} 
& & \int_{-1}^1 (1-x)^{\alpha} (1+x)^{\beta}
P_k^{(\alpha,\beta)} (x)P_l^{(\alpha,\beta)} (x) \; dx=
\nonumber \\ \nonumber
& & \frac{2^{\alpha+\beta+1}}{2l+\alpha+\beta+1}
\frac{\Gamma(l+\alpha+1)\Gamma(l+\beta+1)}{\Gamma(l+\alpha+\beta+1)l!}
\delta_{kl}
\end{eqnarray}
for $\mathrm{Re}(\alpha)>-1$ and $\mathrm{Re}(\beta)>-1$.

The $k$-th derivative of $P_l^{(\alpha,\beta)} (x)$, $k\in\mN_0$, is
\begin{eqnarray}
\frac{d^k}{d x^k}
P_l^{(\alpha,\beta)} (x) =
\frac{\Gamma (\alpha+\beta+l+1+k)}{2^k \Gamma (\alpha+\beta+l+1)}
P_{l-k}^{(\alpha+k, \beta+k)} (x) .
\end{eqnarray}
The Jacobi polynomials $P_l^{(\alpha,\beta)}(x)$ are the polynomial 
solutions of the hypergeometric differential equation

\begin{eqnarray}
\label{eqn:9.3}
\left((1-x^2)\frac{d^2}{dx^2}
 + ( \beta-\alpha - (\alpha + \beta + 2)x )\frac{d}{dx}+
l(l+\alpha+\beta+1)\right) P_l^{(\alpha,\beta)} (x) = 0,\nonumber \\
\end{eqnarray}
and specialize for $\alpha =\beta$ to the Gegenbauer
polynomials fulfilling recurrence relation
\begin{eqnarray}
C_{l}^\alpha(x) = \frac{1}{l}\left(2x(l+\alpha-1)C_{l-1}^\alpha(x) -
(l+2\alpha-2)C_{l-2}^\alpha(x)\right)
\end{eqnarray}
with $C_0^\alpha(x)  = 1, C_1^\alpha(x)  = 2 \alpha x$.
The Gegenbauer polynomials are solutions of the Gegenbauer differential equation
\begin{equation}
\label{eqn:Gdiff}
\left((1-x^{2})\frac{d^2}{dx^2}-(2\alpha+1)x\frac{d}{dx}
+l(l+2\alpha)\right) C_{l}^\alpha(x)=0,
\end{equation}
and are represented by finite hypergeometric series
$$
C_l^{\alpha}(z)=\frac{(2\alpha)_{{l}}}{l!}
\,_2F_1\left(-l,2\alpha+l;\alpha+\frac{1}{2};\frac{1-z}{2}\right).
$$
More explicitly,
\begin{equation}
\label{eqn:Cpoly}
C_l^{\alpha}(z)=\sum_{k=0}^{[ l/2 ]}
(-1)^k\frac{\Gamma(l-k+\alpha)}{\Gamma(\alpha)k!(l-2k)!}(2z)^{l-2k},
\end{equation}
and their relation to the Jacobi polynomials is
\begin{eqnarray}
C_l^{\alpha}(x) =
\frac{(2\alpha)_{l}}{(\alpha+\frac{1}{2})_{l}}P_l^{(\alpha-1/2,\alpha-1/2)}(x).
\end{eqnarray}


\subsection*{Acknowledgments}

P.Pand\v zi\' c  acknowledges the financial support  from the grant no. 4176 by the Croatian Science Foundation.
P. Somberg acknowledges the financial support from the grant GA P201/12/G028.


\vspace{0.5cm}

\flushleft Pavle Pand\v zi\' c\\
Department of Mathematics, University of Zagreb \\
Bijeni\v cka 30, 10000 Zagreb, Croatia \\
E-mail: pandzic@gmail.com 

\vspace{0.2cm}
\flushleft Petr Somberg\\
Mathematical Institute of Charles University\\
Sokolovsk\'a 83, Praha 8 - Karl\'{\i}n, Czech Republic \\
E-mail: somberg@karlin.mff.cuni.cz

\end{document}